\documentclass[12pt,leqno]{article}
\usepackage{amsfonts}
\usepackage{amsmath, amsthm, amsfonts, amssymb, color}
\usepackage{mathrsfs}
\usepackage{color}
\theoremstyle{definition}
\newtheorem{defn}{Definition}[section]
\newcommand{\scr}[1]{\mathscr #1}
\definecolor{wco}{rgb}{0.5,0.2,0.3}

\numberwithin{equation}{section} \theoremstyle{remark}

\newcommand{\ua}{\uparrow}

\title{{\bf   Exponential Convergence  of Non-Linear Monotone SPDEs }\footnote{Supported in
 part by  NNSFC(11131003, 11431014), the 985 project and the Laboratory of Mathematical and  Complex Systems.} }
\author{
{\bf     Feng-Yu Wang  }\\
\footnotesize{ School of Mathematical Sciences,
Beijing Normal
University, Beijing 100875, China}\\
 \footnotesize{ Department of Mathematics,
Swansea University, Singleton Park, SA2 8PP, United Kingdom}\\
\footnotesize{  wangfy@bnu.edu.cn, F.-Y.Wang@swansea.ac.uk}}
\begin{document}
\allowdisplaybreaks
\def\R{\mathbb R}  \def\ff{\frac} \def\ss{\sqrt} \def\B{\mathbf
B}
\def\N{\mathbb N} \def\kk{\kappa} \def\m{{\bf m}}
\def\ee{\varepsilon}\def\ddd{D^*}
\def\dd{\delta} \def\DD{\Delta} \def\vv{\varepsilon} \def\rr{\rho}
\def\<{\langle} \def\>{\rangle} \def\GG{\Gamma} \def\gg{\gamma}
  \def\nn{\nabla} \def\pp{\partial} \def\E{\mathbb E}
\def\d{\text{\rm{d}}} \def\bb{\beta} \def\aa{\alpha} \def\D{\scr D}
  \def\si{\sigma} \def\ess{\text{\rm{ess}}}
\def\beg{\begin} \def\beq{\begin{equation}}  \def\F{\scr F}
\def\Ric{\text{\rm{Ric}}} \def\Hess{\text{\rm{Hess}}}
\def\e{\text{\rm{e}}} \def\ua{\underline a} \def\OO{\Omega}  \def\oo{\omega}
 \def\tt{\tilde} \def\Ric{\text{\rm{Ric}}}
\def\cut{\text{\rm{cut}}} \def\P{\mathbb P} \def\ifn{I_n(f^{\bigotimes n})}
\def\C{\scr C}      \def\aaa{\mathbf{r}}     \def\r{r}
\def\gap{\text{\rm{gap}}} \def\prr{\pi_{{\bf m},\varrho}}  \def\r{\mathbf r}
\def\Z{\mathbb Z} \def\vrr{\varrho} \def\ll{\lambda}
\def\L{\scr L}\def\Tt{\tt} \def\TT{\tt}\def\II{\mathbb I}
\def\i{{\rm in}}\def\Sect{{\rm Sect}}  \def\H{\mathbb H}
\def\M{\scr M}\def\Q{\mathbb Q} \def\texto{\text{o}} \def\LL{\Lambda}
\def\Rank{{\rm Rank}} \def\B{\scr B} \def\i{{\rm i}} \def\HR{\hat{\R}^d}
\def\to{\rightarrow}\def\l{\ell}\def\iint{\int}
\def\EE{\scr E}\def\no{\nonumber}
\def\A{\scr A}\def\V{\mathbb V}
\def\BB{\scr B}\def\Ent{{\rm Ent}}

\maketitle

\begin{abstract}  For a Markov semigroup $P_t$ with invariant probability measure $\mu$, a constant $\ll>0$ is called a lower bound of the ultra-exponential convergence rate of $P_t$ to $\mu$, if there exists a constant $C\in (0,\infty)$ such that
$$ \sup_{\mu(f^2)\le 1}\|P_tf-\mu(f)\|_\infty \le C \e^{-\ll t},\ \ t\ge 1.$$     By using the coupling by change of measure in the line of [F.-Y. Wang,   Ann. Probab. 35(2007), 1333--1350], explicit lower bounds of the   ultra-exponential convergence rate are derived  for  a class of non-linear monotone stochastic partial differential equations. The main result is illustrated by the stochastic porous medium equation and the stochastic $p$-Laplace equation respectively.
Finally, the $V$-uniformly  exponential convergence is investigated for  stochastic fast-diffusion equations.
\end{abstract} \noindent
 \vskip 2cm

\section{Introduction}

It is well known that the solution to the porous medium equation
\beq\label{10}\d X_t= \DD X_t^r\,\d t\end{equation} decays at the algebraic rate $t^{-\ff1{r-1}}$ as $t\to\infty$, where $\DD$ is the Dirichlet Laplacian on a bounded  domain in $\R^d$, $r>1$ is a constant and $X_t^r:=
|X_t|^{r-1}X_t$, see \cite{AP}. This type algebraic convergence has been extended in \cite{DRRW} to stochastic generalized porous media equations. When $r\in (0,1)$,   \eqref{10}  is called the fast-diffusion equation.

Consider, for instance, the stochastic porous medium equation
$$\d X_t= \DD X_t^r\d t+ \d W_t,$$ where $\DD$ is the Dirichlet Laplacian on $(0,l)$ for some $l>0$, and $W_t$ is the cylindrical Brownian motion on $L^2(\m),$ where $\m$ is the normalized Lebesgue measure on $(0,l)$. By \cite[Theorem 1.3]{DRRW}, for any
$x\in \H:= H^{-1}$ (the duality of the Sobolev space w.r.t. $L^2(\m)$, see Section 3), the equation has a unique solution   starting at $x$, and the associated Markov semigroup $P_t$ has a unique invariant probability measure $\mu$ such that
\beq\label{1.0}\|P_tf-\mu(f)\|_{\infty}\le C\mathcal L(f) t^{-\ff1 {r-1}},\ \ t>0\end{equation} holds for some constant $C>0$ and all Lipschitz continuous function $f$, where $\mathcal L(f)$ is the Lipschitz constant of $f$,   $\|f\|_\infty:=\sup_{x\in \H}|f(x)|$
and $\mu(f):=\int_\H f\d\mu$.

On the other hand, by using the dimension-free Harnack inequality and a result due to \cite{GM}, the uniform exponential convergence
$$\|P_t f-\mu(f)\|_\infty \le C\e^{-\ll t}\|f\|_\infty,\ \ t\ge 0, f\in L^2(\mu)$$ is proved in \cite{L} for some constants $C,\ll>0$. Since, according to \cite[Theorem 1.2(4)]{W07} (see also \cite[Theorem 1.4(iv)]{L})  $P_t$ is ultrabounded, i.e. $\|P_t\|_{L^2(\mu)\to L^\infty(\mu)}<\infty$ for $t>0$, this implies the ultra-exponential convergence
\beq\label{1.1} \|P_t f-\mu(f)\|_{\infty}^2 \le C\{\mu(f^2) -\mu(f)^2\} \e^{-\ll t},\ \ t\ge 1, f\in L^2(\mu)\end{equation}   for some constant $C,\ll>0$.   To see that \eqref{1.1} improves \eqref{1.0}  for large time, we note that
\beg{equation*}\beg{split} \mu(f^2)-\mu(f)^2&=\ff 1 2\int_{\H\times \H} |f(x)-f(y)|^2\mu(\d x)\mu(\d y) \\
&\le \mathcal L(f)^2\int_{\H\times \H}|x-y|_\H^2\mu(\d x)\mu(\d y)=:C' \mathcal L(f)^2\end{split}\end{equation*}
with constant $C'\in (0,\infty)$ since $\mu(\|\cdot\|_\H^2)<\infty$, see for instance  \cite[Theorem 1.3]{DRRW}.

 However, explicit estimates on the ultra-exponential convergence rate  $\ll$ is not yet available. We note that in \cite{GM2} an lower bound estimate of exponential convergence rate is presented for a class of semi-linear SPDEs (stochastic partial differential equations). But the main result   in \cite{GM2} does not apply to the present non-linear model, since  both  \cite[Hypothesis 2.4(a)]{GM2} (i.e. $F$ is a Lipschitz map from $\H$ to $\H$) and \cite[Hypothesis 2.4(b)]{GM2} (i.e. ${\rm Im}(F)\subset L^2(\m)$) are not satisfied for the present $F(x):= \DD x^r$, which is not a well defined map from $\H$ to $\H$.

In this paper, we aim  to present explicit lower bound estimates for the ultra-exponential convergence rate $\ll$ in \eqref{1.1}.
In the next section, we prove a general result for a class of non-linear SPDEs considered in \cite{L}. The main tool in the study is the coupling by change of measure constructed in \cite{W07} (see also \cite{L}). A general theory on this kind of couplings and applications has been addressed in the recent monograph \cite{Wbook}. The main result is applied to the stochastic porous medium equation and the stochastic $p$-Laplace equation in Section 3 and Section 4  respectively. Finally, in Section 5 we investigate the exponential convergence for stochastic fast-diffusion equations.

\section{A general result}

Let $\V\subset \H\subset \V^*$ be a Gelfand triple, i.e. $(\H,\<\cdot,\cdot\>_\H, |\cdot|_\H)$ is a separable Hilbert space, $\V$ is a reflexive Banach space continuously and densely embedded into $\H$, and $\V^*$ is the duality of $\V$ with respect to $\H$.  Let $_{\V^*}\<\cdot,\cdot\>_\V$ be the dualization between   $\V$ and   $\V^*$. We have $_{\V^*}\<u,v\>_\V=\<u,v\>_\H$ for $u\in\H$ and $v\in\V.$

 Let $W=(W_t)_{t\ge 0}$ be a cylindrical Brownian motion on a (possibly different) Hilbert space $(E, \<\cdot,\cdot\>_E,|\cdot|_E)$, i.e. $W_t:= \sum_{i=1}^\infty B_t^i e_i$ for an orthonormal basis $\{e_i\}_{i\ge 1}$ of $E$ and a sequence of independent one-dimensional Brownian motions $\{B_t^i\}_{i\ge 1}$ on a complete filtered probability space $(\OO,\F,\{\F_t\}_{t\ge 0},\P).$
 Consider the following stochastic equation:
 \beq\label{2.1} \d X_t= b(X_t)\d t + Q\d W_t,\end{equation}
 where
 $b:  \V \to \V^*$ is measurable and $Q\in \L_{HS}(E;\H),$ the space of all Hilbert-Schmidt linear operators from $E$ to $\H$, such that   the following assumptions hold for some constants $r>1, C_1\ge 0$  and  $C_2>0$:
 \beg{enumerate} \item[{\bf (A1)}] (Hemicontinuity) For any   $v_1,v_2,v\in\V$,\,
 $\R\ni s\mapsto {_{\V^*}\<}b(v_1+s v_2), v\>_\V$ is continuous.
 \item[{\bf (A2)}] (Monotonicity) For any $ v_1,v_2\in\V$,\,
 $_{\V^*}\<b(v_1)-b(v_2), v_1-v_2\>_\V \le C_1|v_1-v_2|_\H^2.$
 \item[{\bf (A3)}] (Coercivity) For any $  v\in\V,$ \,
 $_{\V^*}\<b(v),v\>_\V \le C_1 -C_2 \|v\|_\V^{r+1}.$
 \item[{\bf (A4)}] (Growth) For any $u,v\in \V,$\,
 $ |_{\V^*}\<b(v),u\>_\V|\le C_1 \big(1+\|v\|_\V^{r+1}+\|u\|_\V^{r+1}\big).$ \end{enumerate}

\begin{defn}\label{D2.1} A continuous $\H$-valued adapted process $X$ is called a   solution to (\ref{2.1}), if
$$ \int_0^T \E\|X_t\|_\V^{r+1}\,\d t  <\infty,\ \ T>0,$$   and $\P$-a.s.
$$X_t=X_0 +\int_0^t b(X_s)\,\d s +\int_0^t Q \d W_s,\ \ t\ge 0,$$ where the Bochner integral $\int_0^t b(X_s)\,\d s$  is defined on $\V^*$.   \end{defn}

According to \cite[Theorems II.2.1, II.2.2]{KR}, for any $x\in\H$, the equation (\ref{2.1}) has a unique solution $X_t^x$ with initial datum $x$; see also \cite[Theorem 2.1]{RRW} for
$${\bf K}:= L^{1+r}([0,T]\times\OO\to \V; \d t\times \P)\cap L^2([0,T]\times \OO\to\H;\d t\times\P).$$
Let $P_t$ be the associated Markov semigroup, i.e.
$$P_t f(x):= \E f(X_t^x),\ \ f\in \B_b(\H), t\ge 0, x\in\H.$$ For any $u\in \H$, let
$$\|u\|_Q= \inf\{|x|_E: x\in E, Qx=u\},$$ where we set $\inf\emptyset =\infty$ by convention.

The study of (\ref{2.1}) with the above type assumptions goes back to \cite{P1,P2} for  non-linear  monotone SPDEs.  Extensions to stochastic equations with $``$local conditions" as well as to   non-monotone stochastic equations have been made  in \cite{LR, LR2,RW08}. As mentioned in the Introduction that in this paper we aim to estimate the ultra-convergence rate of $P_t$. The following is the main result of the paper.

\beg{thm}\label{T2.1} Assume that ${\rm Ker}Q=\{0\}.$ If for some   constant  $\theta\in [2,\infty)\cap (r-1,\infty)$ there exist $\eta,\dd\in (0,\infty)$  such that
\beq\label{2.2} 2\,_{\V^*}\<b(u)-b(v), u-v\>_{\V}\le -\max\big\{\eta\|u-v\|_{Q}^\theta |u-v|^{r+1-\theta}_\H,\ \dd|u-v|^{1+r}_\H\big\}\end{equation} holds for all $u,v\in\V,$ then $P_t$ has a unique invariant probability measure $\mu$ and
$\eqref{1.1}$ holds for some constant $C>0$ and
\beq\label{FY}\ll:= \sup_{t>0}\ff 1 t \log \ff 2 {\exp\big[\aa t^{-\ff{r+1}{r-1}}\big]-1} >0,\end{equation} where
$$\aa:= \Big(\ff{\theta r+\theta}{r-1}\Big)^{\ff{r+1}{r-1}}\ff{2+\theta}{(\theta+1-r)^{\ff{2r}{r-1}}\dd^{\ff{2(\theta+1-r)}{\theta(r-1)}}\eta^{\ff 2\theta}}\in (0,\infty).$$
Moreover,
\beg{equation*}\beg{split} \ll &\ge \ff{(\theta+1-r)^{\ff{2r}{r+1}}\dd^{\ff{2(\theta+1-r)}{\theta(r+1)}}\eta^{\ff{2(r-1)}{\theta(r+1)}}}{\theta (2+\theta)^{\ff{r-1}{r+1}}}
 \Big\{\log \big(1+2\e^{-\ff{1+r}{r-1}}\big)\Big\}^{\ff{r-1}{r+1}} \\
&\ge \ff{(\theta+1-r)^{\ff{2r}{r+1}}\dd^{\ff{2(\theta+1-r)}{\theta(r+1)}}\eta^{\ff{2(r-1)}{\theta(r+1)}}}{\e \theta (2+\theta)^{\ff{r-1}{r+1}}}.\end{split}\end{equation*}
\end{thm}

\beg{proof} By  \cite[Theorem 1.4]{L} with $\aa=1+r$ (see also \cite[Corollary 2.2.4]{Wbook} with $\aa=r$), $P_t$ has a unique invariant probability measure $\mu$ of full support on $\H$. Moreover,   $P_t$ is strong Feller (i.e. $P_t \B_b(\H)\subset C_b(\H), t>0$) and ultra-bounded (i.e. $\|P_t\|_{L^2(\mu)\to L^\infty(\mu)}<\infty, t>0$) with
\beq\label{UL} \|P_t\|_{L^2(\mu)\to L^\infty(\mu)}\le \exp\Big[c+ct^{-\ff{r+1}{r-1}}\Big],\ \ t>0\end{equation} holding for some constant $c>0$. So,
$$\|P_t f\|_\infty^2=\|P_t f\|_{L^\infty(\mu)}^2 \le \|P_1\|_{L^2(\mu)\to L^\infty(\mu)}^2 \mu((P_{t-1}f)^2),\ \ t\ge 1. $$
Therefore, it suffices to prove
\beq\label{1.1'} \mu((P_t f)^2)\le C\mu(f^2) \e^{-\ll t},\ \ t\ge 0,f\in L^2(\mu), \mu(f)=0\end{equation} for some constant $C\in (0,\infty)$ and the desired constant $\ll$,  and to verify the claimed lower bounds of $\ll$.
 We shall complete the proof by four steps.

\

Step 1.  We first construct a coupling by change of measure using the idea of \cite{W07}. For fixed $T>0$ and $x,y\in \H$, let $X_t=X_t^x$ solve \eqref{2.1} for $X_0=x$, and let $Y_t$ solve the equation
\beq\label{Y} \d Y_t= \Big\{b(Y_t)+\ff{\bb (X_t-Y_t)}{|X_t-Y_t|^\vv_\H}\Big\}\d t +Q\d W_t,\ \ Y_0=y,\end{equation} where
\beq\label{BE} \vv:= \ff{\theta+1-r}{2+\theta}\in (0,1),\ \ \ \bb:= \ff{|x-y|^\vv_\H}{\vv T}\ge 0,\ \ \text{and}\ \ff{X_t-Y_t}{|X_t-Y_t|_\H^\vv}:=0\ \text{if}\ X_t=Y_t.\end{equation}
As shown in \cite[Theorem A.2]{W07} (see also \cite{L}) that   the equation \eqref{Y} has a unique solution such that $X_t=Y_t$ for $t\ge\tau$, where $$\tau:=\inf\{t\ge 0: X_t=Y_t\}$$ is the coupling time.
By \eqref{2.2} we have
\beq\label{W*0} \d |X_t-Y_t|_\H^2\le -\big\{\eta\|X_t-Y_t\|_Q^\theta|X_t-Y_t|_\H^{r+1-\theta}+2\bb |X_t-Y_t|_\H^{2-\vv}\big\}\d t,\ \ t<\tau.\end{equation}
Since $\eta>0$, this implies
$$\d |X_t-Y_t|_\H^\vv =\ff\vv 2 (|X_t-Y_t|_\H^2)^{\ff{\vv-2}2}\d|X_t-Y_t|_\H^2\le -\vv\bb\d t,\ \ t<\tau.$$ Thus, if $T<\tau$ then
$$|X_T-Y_T|_\H^\vv\le |x-y|_\H^\vv-\vv \bb T=0,$$ which is a contradiction since by definition it implies $\tau\le T$.    Therefore, we have $\tau\le T$, so that $X_T=Y_T$.

\

Step 2.   By \eqref{W*0}, $\bb\ge 0$ and noting that $2(\vv-1)+r+1-\theta=-\vv\theta$, we have
$$\d|X_t-Y_t|_\H^{2\vv}=\vv |X_t-Y_t|_\H^{2(\vv-1)}\d |X_t-Y_t|_\H^2 \le-\vv\eta \ff{\|X_t-Y_t\|_Q^\theta}{|X_t-Y_t|_\H^{\vv\theta}}\,\d t,\ \ \ t<\tau.$$ Then
\beq\label{EW}\eta\int_0^\tau \ff{\|X_t-Y_t\|_Q^\theta}{|X_t-Y_t|_\H^{\vv\theta}}\,\d t\le \ff{|x-y|_\H^{2\vv}}\vv =\ff{2+\theta}{\theta+1-r}|x-y|_\H^{\ff{2(\theta+1-r)}{2+\theta}}.\end{equation}  Combining this with $T\le \tau$ and $\theta\ge 2$, and using the Jensen inequality, we see that
$$\zeta_t:= \ff{\bb Q^{-1}(X_t-Y_t)}{|X_t-Y_t|_\H^\vv} $$ is well defined in $L^2([0,T]\to E;\d t)$. Moreover, since $\theta\ge 2$, by \eqref{EW} and the H\"older inequality we obtain
\beq\label{W*}\beg{split} &\int_0^T|\zeta_t|_E^2\d t =\int_0^\tau \ff{\bb^2\|X_t-Y_t\|_Q^2}{|X_t-Y_t|_\H^{2\vv}}\,\d t\\
 &\le \bigg(\eta\int_0^\tau \ff{\|X_t-Y_t\|_Q^\theta}{|X_t-Y_t|_\H^{\vv\theta}}\,\d t\bigg)^{\ff 2\theta} \bigg(\int_0^T\ff{\bb^{\ff {2\theta}{\theta-2}}}{\eta^{\ff 2{\theta-2}}}\,\d t\bigg)^{\ff{\theta-2}\theta}\\
&\le  \ff{T^{\ff{\theta-2}\theta}\bb^2}{\eta^{\ff 2\theta}} \Big(\ff{2+\theta}{\theta+1-r}\Big)^{\ff 2\theta} |x-y|_\H^{\ff{4(\theta+1-r)}{\theta(2+\theta)}}\\
&=  \ff{(2+\theta)^{\ff{2(\theta+1)}{\theta}}|x-y|_\H^{\ff{2(\theta+1-r)}\theta}}{(\theta+1-r)^{\ff{2(\theta+1)}\theta} T^{\ff{\theta+2}\theta}\eta^{\ff 2\theta}}.\end{split}\end{equation}
 Then by the Girsanov theorem,
$$R:=\exp\bigg[-\int_0^T\<\zeta_t,\d W_t\>_E-\ff 1 2 \int_0^T |\zeta_t|_E^2\d t\bigg]$$ is a well defined probability density of $\P$, and the process
$$\tt W_t:= W_t+\int_0^t \zeta_s\d s,\ \ t\in [0,T]$$ is a cylindrical Brownian motion on $E$ under the weighted probability measure $\d\Q:=R\d\P.$ Now, rewrite \eqref{Y} by
$$\d Y_t = b(Y_t)\d t+Q\d \tt W_t,\ \ \ Y_0=y.$$ From  the weak uniqueness of the solution to \eqref{2.1} and  $X_T=Y_T$, we conclude that
$$P_Tf(y)= \E_\Q f(Y_T)= \E[Rf(Y_T)]= \E[Rf(X_T)].$$ This together with $P_Tf(x)= \E f(X_T)$ yields that
\beq\label{*1} |P_Tf(x)-P_Tf(y)|^2 = \big|\E[f(X_T)(1-R)]\big|^2 \le (P_T f^2(x))(\E R^2-1).\end{equation}

\

Step 3.  By \eqref{W*} we have
\beq\label{W*2} \beg{split} &\E R^2 = \E \exp\bigg[-2\int_0^T\<\zeta_t,\d W_t\>_E-  \int_0^T |\zeta_t|_E^2\d t\bigg]\\
&\le \exp\bigg[\ff{(2+\theta)^{\ff{2(\theta+1)}\theta}|x-y|_\H^{\ff{2(\theta+1-r)}\theta}}{(\theta+1-r)^{\ff{2(\theta+1)}\theta}T^{\ff{\theta+2}\theta}\eta^{\ff 2\theta}}\bigg]\E  \e^{-2\int_0^T\<\zeta_t,\d W_t\>_E-  2\int_0^T |\zeta_t|_E^2\d t}\\
&= \exp\bigg[\ff{(2+\theta)^{\ff{2(\theta+1)}\theta}|x-y|_\H^{\ff{2(\theta+1-r)}\theta}}{(\theta+1-r)^{\ff{2(\theta+1)}\theta}T^{\ff{\theta+2}\theta}\eta^{\ff 2\theta}}\bigg].\end{split}\end{equation} Moreover,   \eqref{2.2} yields that
 $$\d |X_t^x-X_t^y|_\H^2 \le -\dd   |X_t^x-X_t^y|_\H^{1+r}\d t,\ \ \ t\ge 0, x,y\in \H,$$ where $X_t^x$ and $X_t^y$ solve the equation \eqref{2.1} starting at $x$ and $y$ respectively.  Thus,
$$ |X_t^x-X_t^y|_\H^2 \le \Big(\ff {\dd t(r-1)} 2 \Big)^{\ff 2 {1-r}},\ \ t>0, x,y\in \H.$$
Substituting this and \eqref{W*2} into \eqref{*1} and using the Markov property,   we arrive at
\beg{equation}\label{ED}\beg{split} &|P_{T+s}f(x)-P_{T+s}f(y)|^2 \le \E |P_T f(X_s^x)-P_T f(X_s^y)|^2 \\
&\le \E\bigg\{P_{T}f^2(X_s^x) \bigg(\exp\bigg[\ff{(2+\theta)^{\ff{2(\theta+1)}\theta}|X_s^x-X_s^y|_\H^{\ff{2(\theta+1-r)}\theta}}{(\theta+1-r)^{\ff{2(\theta+1)}\theta}T^{\ff{\theta+2}\theta}\eta^{\ff 2\theta}}\bigg]-1\bigg)\bigg\}\\
&\le (P_{T+s}f^2(x)) \bigg(\exp\Big[\ff{C_0}{s^{\ff{2(\theta+1-r)}{\theta(r-1)}} T^{\ff{\theta+2}\theta}}\Big]-1\bigg),\ \ T,s>0,\end{split}\end{equation}
where
\beq\label{*3} C_0:= \eta^{-\ff 2\theta}  \Big(\ff{2+\theta}{\theta +1-r}\Big)^{\ff{2(\theta+1)}{\theta}} \Big(\ff 2 {\dd(r-1)}\Big)^{\ff{2(\theta+1-r)}{\theta (r-1)}}.\end{equation} For fixed $t>0$, by taking $s\in (0,t)$ and $T=t-s$ in \eqref{ED}  we obtain
\beq\label{0*} \beg{split} \mu((P_tf)^2) &= \ff 1 2 \int_{\H\times \H} |P_tf(x)-P_tf(y)|^2\mu(\d x)\mu(\d y)\\
&\le \ff{\mu(f^2)}2 \inf_{s\in (0,t)} \bigg\{\exp\Big[\ff{C_0}{s^{\ff{2(\theta+1-r)}{\theta(r-1)}}(t-s)^{\ff{2+\theta}{\theta}}}\Big]-1\bigg\},\ \ t>0,\ \mu(f)=0.\end{split}\end{equation}

Step4.  To calculate the inf in \eqref{0*}, let $$\aa_1= \ff{2(\theta+1-r)}{\theta(r-1)}, \ \aa_2= \ff{2+\theta}{\theta}.$$
We have  $\aa_1+\aa_2= \ff{r+1}{r-1},$ and by \eqref{*3},
\beg{equation*}\beg{split} &\inf_{s\in (0,t)} \ff{C_0}{s^{\aa_1}(t-s)^{\aa_2}} = \ff{C_0 (\aa_1+\aa_2)^{\aa_1+\aa_2}}{t^{\aa_1+\aa_2} \aa_1^{\aa_1}\aa_2^{\aa_2}}\\
&= \ff{C_0 }{t^{\ff{r+1}{r-1}}} \Big(\ff{r+1}{r-1}\Big)^{\ff{r+1}{r-1}} \Big(\ff{\theta(r-1)}{2(\theta+1-r)}\Big)^{\ff{2(\theta+1-r)}{\theta(r-1)}}
\Big(\ff{\theta}{2+\theta}\Big)^{\ff{2+\theta}{\theta}}\\
&= \ff{\aa}{t^{\ff{r+1}{r-1}}},\ \ \ \ t>0.\end{split}\end{equation*}
Then it follows from \eqref{0*} that
\beq\label{1*} \mu((P_tf)^2)\le \ff {\mu(f^2)}2 \bigg(\exp\Big[\aa t^{-\ff{r+1}{r-1}}\Big]-1\bigg),\ \ \ t>0, \mu(f)=0.\end{equation}
Obviously,   there exists $t_0\in (0,\infty)$ such that
$$0<\ff 1 {t_0} \log \ff 2{\exp\big[\aa t_0^{-\ff{r+1}{r-1}}\big]-1}=\sup_{t>0}  \ff 1 {t} \log \ff 2{\exp\big[\aa t^{-\ff{r+1}{r-1}}\big]-1}=:\ll.$$ So, \eqref{1*} yields that
$$\mu\big((P_{t_0}f)^2\big) \le \mu(f^2)\e^{-\ll t_0},\ \ \mu(f)=0.$$
Letting
$i(t)=\sup\{n\in \Z_+: n\le \ff t {t_0}\}$ be the integer part of $\ff t{t_0},$ combining this with    the semigroup property and the $L^2(\mu)$-contraction of $P_t$, we obtain
$$\mu((P_tf)^2)\le \mu((P_{i(t) t_0}f)^2)  \le \mu(f^2)\e^{-\ll t_0 i(t)} \le \mu(f^2) \e^{-\ll (t- t_0)},\ \ t\ge 0,\mu(f)=0.$$ Thus, \eqref{1.1'} holds for $C:= \e^{\ll t_0}.$

Finally, to derive the desired explicit lower bounds of $\ll$, we take
$$t= \Big(\ff {\aa}{\log (1+2\exp[-\ff{r+1}{r-1}])}\Big)^{\ff{r-1}{r+1}}.$$
Then
\beg{equation*}\beg{split} \ll&\ge \ff 1{t} \log \ff 2{\exp\big[\aa t^{-\ff{r+1}{r-1}}\big]-1} =\ff{(r+1)\{\log  (1+2\exp\big[-\ff{r+1}{r-1}\big])\}^{\ff{r-1}{r+1}}}{(r-1)\aa^{\ff{r-1}{r+1}}}\\
&=  \ff{(\theta+1-r)^{\ff{2r}{r+1}}\dd^{\ff{2(\theta+1-r)}{\theta(r+1)}}\eta^{\ff{2(r-1)}{\theta(r+1)}}}{\theta (2+\theta)^{\ff{r-1}{r+1}}}
 \Big\{\log \big(1+2\e^{-\ff{1+r}{r-1}}\big)\Big\}^{\ff{r-1}{r+1}} \\
&\ge \ff{(\theta+1-r)^{\ff{2r}{r+1}}\dd^{\ff{2(\theta+1-r)}{\theta(r+1)}}\eta^{\ff{2(r-1)}{\theta(r+1)}}}{\e \theta (2+\theta)^{\ff{r-1}{r+1}}},\end{split}\end{equation*}  where the last step is due to the fact that
$$\inf_{s\ge 1} \{\log (1+ 2 \e^{-s})\}^{\ff 1 s} =\lim_{s\to\infty} \{\log (1+ 2 \e^{-s})\}^{\ff 1 s}=\e^{-1}.$$
\end{proof}

To conclude this section, we indicate that   $P_t$ is ultra-exponential convergent provided
$$2\,_{\V^*}\<b(u)-b(v), u-v\>_{\V}\le \gg |u-v|_\H^2 -\max\big\{\eta\|u-v\|_{Q}^\theta |u-v|^{r+1-\theta}_\H,\ \dd|u-v|^{1+r}_\H\big\}$$
holds for some constant $\gg,\eta>0$, which is weaker than \eqref{2.2}. This can be proved as   in \cite[proof of Theorem 1.5]{L}  using the Harnack inequality in \cite[Theorem 2.2.1]{Wbook} and the ultraboundedness of $P_t$.  When $\gg>0$ is small enough, with the coupling constructed in  the proof of Theorem 2.2.1 in \cite{Wbook}, we may  derive explicit lower bounds of the convergence rate $\ll$ using the argument   in the proof of Theorem \ref{T2.1}. As in this case the resulting estimates will be rather complicated, in Theorem \ref{T2.1}   we only consider the case that $\gg=0.$ However,   to derive explicit lower bounds of $\ll$ for any $\gg>0$,  new techniques are required.

\section{Stochastic  porous medium equation}

Let $\DD$ be the Dirichlet Laplacian on the interval $(0,l)$ for some $l>0$, and let $\si>0, r>1$ be two constants. Let $W_t$ be the cylindrical Brownian motion on $L^2(\m)$, where $\m(\d x):=l^{-1}\d x$ is the normalized Lebesgue measure on $(0,l).$ Consider the following stochastic porous medium equation
\beq\label{E3} \d X_t= \DD X_t^r \d t +\si \d W_t.\end{equation}

We first verify assumptions {\bf (A1)}-{\bf (A4)} for an appropriate choice of $(\H,\V)$.
It is well known that the spectrum of $-\DD$ consists of simple eigenvalues   $\{\ll_k:= \ff{\pi^2k^2}{l^2}\}_{k\ge 1}$. Let $\{e_k\}_{k\ge 1}$ be the corresponding eigenbasis.  Then $Q:= \si I$ is Hilbert-Schmidt from $L^2(\m)$ to $\H:=H^{-1}$,
 the completion of $L^2(\m)$ under the inner product
$$\<x,y\>_\H:=\sum_{i=1}^\infty \ff 1 {\ll_i} \<x,e_i\>\<y, e_i\>.$$ Let $\V= L^{1+r}(\m).$ Then $b(v) := \DD v^r$ extends to a unique map from $\V$ to $\V^*$ with
$$_{\V^*}\<b(v),u\>_\V= -\int_0^l v^r u\d\m,\ \ u,v\in \V.$$ This implies {\bf (A3)} and {\bf (A4)} for $C_1=C_2=1.$ Moreover,
  for any $v_1,v_2,v\in \V$,
 $$_{\V^*}\<b(v_1+sv_2), v\>_\V= -\int_0^l (v_1+s v_2)^r v\d\m$$ is continuous in $s\in\R$; that is, {\bf (A1)} holds.
 Finally, we have (see the proof of Proposition \ref{P3.1} below)
\beq\label{3.1} (s^r-t^r)(s-t)\ge 2^{1-r} |s-t|^{1+r},\ \ \ s,t\in\R.\end{equation} Then
 $$_{\V^*}\<b(v_1)-b(v_2),v_1-v_2\>_\V\le -2^{1-r}\|v_1-v_2\|_\V^{1+r},\ \ v_1,v_2\in\V,$$ so that {\bf (A2)} holds for any positive constant $C_1.$
 Therefore,
for any initial point $x\in \H$ the equation \eqref{3.1} has a unique solution starting at $x$. Let $P_t$ be  the associated Markov semigroup.

\beg{prp}\label{P3.1} For the equation $\eqref{E3}$, $P_t$ has a unique invariant probability measure $\mu$ such that   $\eqref{1.1}$ holds for some constant  $C>0$ and $\ll$ defined in $\eqref{FY}$ for
$$\aa:= \ff{l^{\ff 4{r-1}}  (3+r)(r+1)^{\ff{2(r+1)}{r-1}}}{(2\pi)^{\ff 4{r-1}}\si^2 (r-1)^{\ff{r+1}{r-1}}}.$$
Moreover,
\beg{equation*}\beg{split} \ll &\ge \ff{(2\pi)^{\ff 4{r+1}}\si^{\ff{2(r-1)}{r+1}}\{\log(1+2\exp[-\ff{r+1}{r-1}])\}^{\ff{r-1}{r+1}}}{(r+1)l^{\ff 4{r+1}}(3+r)^{\ff{r-1}{r+1}}}\\
&\ge \ff{(2\pi)^{\ff 4{r+1}}\si^{\ff{2(r-1)}{r+1}}}{\e (r+1)l^{\ff 4{r+1}} (3+r)^{\ff{r-1}{r+1}}}.\end{split}\end{equation*} \end{prp}

\beg{proof} We first prove \eqref{3.1}. Obviously, we may assume that $s\lor t\ge 0$, otherwise simply use $-s,-t$ to replace $s,t$ respectively. Moreover, since the positions of $s$ and $t$ are symmetric, we may assume further that $s>t$ (hence, $s\ge 0$). Assuming $s>t$ and $s\ge 0$, we prove \eqref{3.1} by considering the following two situations respectively.

(i) $s>t\ge 0$. Since $0\le s\mapsto s^r$ is convex, we have
$$\ff{\d}{\d s} \ff{s^r-t^r}{(s-t)^r} =\ff{rt}{(s-t)^{r+1}}\big(t^{r-1}-s^{r-1}\big)\le 0.$$ So,
$$\inf_{s>t} \ff{s^r-t^r}{(s-t)^r} = \lim_{s\to\infty} \ff{s^r-t^r}{(s-t)^r}=1,\ \ t\ge 0.$$ Then \eqref{3.1} holds since $2^{1-r}\le1$.

(ii) $s\ge 0>t.$ By the Jensen inequality we have
$$s^r-t^r = 2\Big(\ff{s^r}2 + \ff{|t|^r} 2\Big) \ge 2\Big(\ff{s+|t|}2\Big)^r = 2^{1-r}(s+|t|)^r= 2^{1-r}(s-t)^r.$$ Thus, \eqref{3.1} holds.

\

Now, let $b(x)=\DD x^r, x\in \V:=L^{r+1}(\m)$. Since $Q=\si I$, we have $\|\cdot\|_Q=\ff 1 \si\|\cdot\|_2.$  Combining this with $\|\cdot\|_{r+1}\ge \|\cdot\|_2$, $\ll_1=\ff{\pi^2}{l^2}$ and the definition of $|\cdot|_\H$,   we obtain
$$\|x\|_{r+1} \ge \|x\|_2= \max\big\{\ss{\ll_1}\,  |x|_\H,\ \si\|x\|_Q\big\}=\max\Big\{\ff{\pi}l |x|_\H,\ \si\|x\|_Q\Big\}.$$   Then, due to \eqref{3.1},   for any $\theta\in (r-1,r+1]\cap [2,r+1],$
\beg{equation*}\beg{split} &2\,_{\V^*}\<b(x)-b(y), x-y\>_\V=- 2\int_0^l(x^r-y^r)(x-y)\d \m \le -2^{2-r} \|x-y\|_{r+1}^{r+1}\\
 &\le -\max\big\{\eta \|x-y\|_Q^\theta|x-y|_\H^{r+1-\theta},\ \dd|x-y|_\H^{r+1}\big\},\ \ x,y\in\V:=L^{r+1}(\m)\end{split}\end{equation*}  holds for $$\eta:= 2^{2-r} \si^\theta \Big(\ff\pi l\Big)^{r+1-\theta},\ \ \dd:= 2^{2-r} \Big(\ff\pi l\Big)^{r+1}.$$Therefore, by Theorem \ref{T2.1}, \eqref{1.1} holds for some constant $C\in (0,\infty)$ and
$$\ll:= \sup_{t>0, \theta\in (r-1,r+1]\cap [2,r+1]} \ff 1 t \log \ff 2{\exp\big[\aa_\theta t^{-\ff{r+1}{r-1}}\big]-1},$$ where
\beg{equation*}\beg{split} \aa_\theta &:=  \ff{4^{\ff{r-2}{r-1}}(2+\theta)}{\si^2(\theta+1-r)^{\ff{2r}{r-1}}}\Big(\ff l\pi\Big)^{\ff 4{r-1}}\Big(\ff{\theta r+\theta}{r-1}\Big)^{\ff{r+1}{r-1}}\\
&= \ff{4^{\ff{r-2}{r-1}}}{\si^2} \Big(\ff l\pi\Big)^{\ff 4{r-1}}\Big(\ff{r+1}{r-1}\Big)^{\ff{r+1}{r-1}} \Big(\ff{\theta}{\theta+1-r}\Big)^{\ff{r+1}{r-1}}\ff{2+\theta}{\theta+1-r}.\end{split}\end{equation*}
Noting that $r\ge 1$ implies $\theta+1-r\le\theta$,  so that $\aa_\theta$ is decreasing in $\theta$, we obtain
$$\inf_{\theta\in (r-1,r+1]\cap [2,r+1]}\aa_\theta= \aa_{r+1} = \ff{l^{\ff 4{r-1}}(3+r)(r+1)^{\ff{2(r+1)}{r-1}}}{(2\pi)^{\ff 4{r-1}}\si^2 (r-1)^{\ff{r+1}{r-1}}}=:\aa.$$
So, \eqref{1.1} holds for some $C\in (0,\infty)$ and the desired $\ll$. Moreover,
  as in the proof of Theorem \ref{T2.1} that the desired lower bound estimates follows by taking in \eqref{FY}
$$t = \Big(\ff {\aa}{\log (1+2\exp[-\ff{r+1}{r-1}])}\Big)^{\ff{r-1}{r+1}}.$$
\end{proof}

To conclude this section, let us recall a corresponding result in the linear case, i.e. $r=1$. Let $R=\si I$ and $T_t=\e^{t\DD}.$
In this case, for any $p>2$ there exist constants $C_p,t_p\in (0,\infty)$ such that
 \beq\label{H} \|P_t-\mu\|_{L^2(\mu)\to L^p(\mu)} \le C_p \e^{-\ll_1 t},\ \ t\ge t_p.\end{equation} To see this, we observe that $\si W_t$ is a Wiener process on $\H$ with variance operator $Q e_i:= \ff{\si^2}{\ll_i}, i\ge 1$. Taking $M=0, R=Q$ and $T_t= \e^{t\DD}$, we see that assumptions in \cite[Coroolary 1.4]{RW} hold for $h_1(t)=\e^{-\ll_1 t/2}$ and $h_2(t)=0,$ so that
 \beq\label{KE} \|P_t-\mu\|_{L^2(\mu)\to L^2(\mu)} \le \e^{-\ll_1 t},\ \ t\ge 0.\end{equation}
 Moreover, according to \cite[Theorem 4 c)]{CG}, $P_t$ is hypercontractive, i.e. for any $p>2$ there exists a constant $t_p>0$ such that $\|P_t\|_{L^2(\mu)\to L^p(\mu)}=1$ holds for $t\ge t_p$. Combining this with \eqref{KE} we prove \eqref{H}. Note that in this linear case $P_t$ is not ultra-bounded, so that
   we do not have the ultra-exponential convergence as in \eqref{1.1}.

A feature in the linear case is that the exponential convergence rate $\ll_1$ is independent of $\si$. Note that for $r>1$ the lower bound estimates of $\ll$ presented in Proposition \ref{P3.1} are increasing to $\infty$ as $\si\uparrow \infty$. But if we let $r\downarrow 1$   in these estimates, the lower bounds of $\ll$ tend  to $\ff{2\ll_1}\e$ (of course, the other constant $C$ will tend to $\infty$ since $P_t$ is not ultracontractive for $r=1$),  which is also independent of $\si$.  This indicates that the power of $\si$ included in the lower bound estimates of $\ll$ presented in Proposition \ref{P3.1} is suitable when $r$ goes down to $1$.

\section{Stochastic $p$-Laplace equation}

Again let   $D=(0,l)$ for some $l>0$ and   $\m$ be  the normalized volume measure.   For $p> 2$, let  $\H_0^{1,p}$ be the closure of $C_0^\infty(D)$ with respect to the norm
$$\|f\|_{1,p}:= \|f\|_{p}  +\|\nn f\|_{p},$$ where, since $D$ is one-dimensional, $\nn f:=f'$.  The $p$-Laplacian on $D$ is defined by
$$\DD_p f= \nn\big(|\nn f|^{p-2}\nn f\big),\ \ \ f\in C^2(D).$$
Consider the SPDE
\beq\label{4.1} \d X_t = \DD_p X_t \d t +Q \d W_t,\end{equation} where $W_t$ is a cylindrical Brownian motion on $L^2(\m)$, and $Q\in\scr L(\H)$ is such that
\beq\label{GE} Qe_i= q_i e_i,\ q_i^2\ge \ff{\si^2}{i^2},\ \ \sum_{i=1}^\infty q_i^2<\infty\end{equation} holds for some constants $\si>0$ and $\{q_i\}_{i\ge 1}\subset \R$, recall that $\{e_i\}_{i\ge 1}$ is the eigenbasis of $-\DD$ with respect to the eigenvalues $\ll_i:=\ff{(i\pi)^2}{l^2}, i\ge 1.$

To apply Theorem \ref{T2.1}, let    $E=\H=L^2(\m), \V=\H_0^{1,p}$ and $r=p-1>1.$  Then $Q\in \scr L_{HS}(E,\H)$ by \eqref{GE}, and $b:=\DD_p$ extends to a unique operator from $\V$ to $\V^*$ with
\beq\label{AD} _{\V^*}\<b(v), u\>_\V:=  -\int_0^l |\nn v|^{p-2}\<\nn v,\nn u\>\d\m,\ \ u,v\in \V.\end{equation}
Thus, {\bf (A1)} and {\bf (A4)} with $C_1=1$ hold.
Next, since for any $f\in C_0^\infty(D)$ we have
\beg{equation*}\beg{split} \int_0^l |f|^p\d\m &=\int_0^l \bigg|\int_0^x f'(s)\d s\bigg|^p\m(\d x) \le \int_0^l x^{p-1}\m(\d x) \int_0^x |f'(s)|^p\d s\\
&\le \bigg(\int_0^l |f'|^p\d\m\bigg)\int_0^l x^{p-1}\d x =\ff{l^p}{p} \int_0^l |f'|^p\d\m,\end{split}\end{equation*} it follows that
$$ \|f\|_{1,p} \le (1+lp^{-\ff 1 p}) \|\nn f \|_p,\ \ f\in \V.$$ From  this and \eqref{AD} it is easy to see that   {\bf (A3)}  holds for
 $C_1=0$ and some $C_2>0$. Moreover, according to the first display on page 767 in \cite{L},  
  \beq\label{DD} 2\,_{\V^*}\<b(u)-b(v), u-v\>_\V\le -C\|\nn(u-v)\|_2^p,\ \ u,v\in \V \end{equation} holds for some constant $C>0.$  Then {\bf (A2)} holds for $C_1=0$.   Therefore,  for any $x\in \H$ the equation \eqref{4.1} has a unique solution starting at $x$. Let $P_t$ be the Markov semigroup associated to \eqref{4.1}.

\beg{prp}\label{P4.1} For the equation $\eqref{4.1}$, $P_t$ has a unique invariant probability measure $\mu$ such that   $\eqref{1.1}$ holds for some constant  $C>0$ and $\ll$ defined in $\eqref{FY}$ for
$$\aa:= \Big(\ff{p^2l^2}{\pi^2(p-2)}\Big)^{\ff p{p-2}} \ff{2+p}{\si^2 2^{\ff{4(p-1)}{p-2}}}. $$
Moreover,
$$ \ll  \ge \ff{\pi^2 2^{\ff{4(p-1)}{p}} \si^{\ff{2(p-2)}p}}{p l^2 (2+p)^{\ff{p-2}p}}\big\{\log\big(1+2\e^{-\ff{p}{p-2}}\big)\big\}^{\ff{p-2}p}
 \ge \ff{\pi^2 2^{\ff{4(p-1)}{p}} \si^{\ff{2(p-2)}p}}{\e p l^2 (2+p)^{\ff{p-2}p}}.$$ \end{prp}

\beg{proof} By the Poincar\'e inequality  we have  \beq\label{D1}\m(|\nn(u-v)|^2) \ge \ff{\pi^2}{l^2} \|u-v\|_2^2=\ff{\pi^2}{l^2} \|u-v\|_\H^2.\end{equation}
Next, by \eqref{GE}
$$\|\nn(u-v)\|_2^2=\ff{\pi^2}{l^2} \sum_{i=1}^\infty i^2\<u-v, e_i\>^2_\H\ge \ff{\pi^2\si^2}{l^2}\sum_{i=1}^\infty \ff 1 {q_i^2}\<u-v,e_i\>_\H^2=\ff{\pi^2\si^2}{l^2}\|u-v\|_Q^2.$$
Combining this with \eqref{DD} and \eqref{D1}, we arrive at
$$2\,_{\V^*}\<b(u)-b(v),u-v\>_{\V} \le -2^{p-1}\Big(\ff \pi l\Big)^{p} \max\big\{\si^p \|u-v\|_Q^p,\ \|u-v\|_\H^p\big\}.$$  This implies  \eqref{2.2} for
$$r=p-1, \ \ \theta=p,\ \ \eta= 2^{p-1}\Big(\ff{\pi\si}l\Big)^p,\ \ \dd= 2^{p-1} \Big(\ff \pi l\Big)^p.$$ Therefore, the proof is finished by Theorem \ref{T2.1}.\end{proof}

\section{Exponential convergence for stochastic fast-diffusion equations}

Consider, for instance,  the equation \eqref{E3} in Section 3 for $r\in (0,1)$, i.e. the stochastic fast-diffusion equation. In this case, we do not have the   ultra-exponential convergence, but we are able to derive a weaker version of exponential convergence by combining the Harnack inequality with a result of \cite{GM}, see \cite{L} for the study of the equation for $r\ge 1.$  To  see the difference between the case of $r\ge 1$ and that of $r\in (0,1)$, we come back to the specific equation \eqref{E3}. When $r\ge 1$ all assumptions in  \cite[Theorem 2.5]{GM} can be easily verified (see the proof of Theorem 1.5 in \cite{L}), but when $r\in (0,1)$ one needs additional conditions (see \eqref{N} below) which exclude the equation \eqref{E3} where $Q:= \si I$ for some $\si>0.$

We would like to mention that  in \cite{Liu11} the ergodicity has been investigated  for equations of type 
$$\d X_t= (\DD X_t^r -\gg \|X_t\|_\H^{q-2}X_t)\d t +\si \d W_t,$$ where $\gg\ge 0$ and $q\ge 2$ are constants. Since $r\in (0,1)$, when $\gg>0$ the term $\gg \|X_t\|_\H^{q-2}X_t$ becomes leading in the study of the convergence rate. For instance, according to \cite[Theorems 1.3-1.4]{Liu11}, 
in this case the solution is exponentially ergodic for $q\ge 2$ and uniformly ergodic for $q>2$. Moreover, algebraic convergence of the semigroup has been proved in \cite{LT11} for a class of equations with weakly  dissipative drifts.  

\

From now on, we let $r\in (0,1)$ and consider the equation \eqref{2.1} such that assumptions {\bf (A1)}-{\bf (A4)} hold. Let $P_t$ be the associated Markov semigroup. We aim to investigate the $V$-uniformly exponential convergence
\beq\label{VE} \|P_t-\mu\|_V:= \sup_{|f|\le V} \Big\|\ff{|P_t f-\mu(f)|}V\Big\|_\infty\le C\e^{-\ll t},\ \ t\ge 0\end{equation} for some constants $C,\ll>0$, where $\mu$ is the invariant probability measure of $P_t$ and $V\ge 1$ is a continuous function on $\H$. Obviously, \eqref{VE} is equivalent to
$$\sup_{|f|\le V} |P_tf(x)-\mu(f)|\le CV(x)\e^{-\ll t},\ \ t\ge 0, x\in\H$$ used in \cite[Definition 2.3]{GM}.

\beg{thm}\label{T5.1} If there exists a non-negative measurable function $h$ on $\V$ such that $\{h\le R\}$ is relatively compact in $\H$ for any $R>0$, and \beq\label{5.1}  _{\V^*}\<b(u),u\>_\V\le \aa-\eta \{h(u)\lor \|u\|_\H\}^{1+r},\ \ u,u\in\H, \end{equation}
\beq\label{5.2}   _{\V^*}\<b(u)-b(v),u-v\>_\V\le -\ff{\eta \|u-v\|_Q^\theta}{|u-v|_\H^{\theta-2}\{h(u)\lor h(v)\}^{1-r}},\ \ u,v\in\V\end{equation} hold for some constants $\aa,\eta>0$ and $\theta\ge \ff 4{1+r}.$ Then $P_t$ has a unique invariant probability measure $\mu$, and for any $\gg>0$, there exist two constants $C,\ll>0$ such that $\eqref{VE}$ holds for $V:=\exp[\gg (1+|\cdot|_\H^2)^{\ff{1-r}2}].$\end{thm}

\beg{proof} By \eqref{5.1} and the It\^o formula, we see that
\beq\label{CM}\ff 1 n\int_0^n \E h(X_t^0)^{1+r}\d t \le \ff{2\aa+\|Q\|_{HS}^2}{2\eta}<\infty,\ \ n\ge 0.\end{equation} Since $h$ has relatively compact level sets in $\H$, this implies that the sequence $\{\ff 1 n\int_0^n \dd_0P_t\d t\}_{n\ge 1}$ is tight and each of its weak limit point gives rise to an invariant probability measure of $P_t$. Now, according to the proof of \cite[Theorem 2.5(1)]{GM}, it suffices to verify
\beg{enumerate} \item[(i)] (Assumption 2.1 in \cite{GM}): $P_t$ is strong Feller (i.e. $P_t \B_b(\H)\subset C_b(\H), t>0$) and $P_t 1_U(x)>0$ holds for any $t>0, x\in \H$ and non-empty open set $U\subset \H$.
\item[(ii)] (Assumption 2.2 in \cite{GM}):  For any $r>0$ there exists $t_0>0$ and a compact subset $K$ of $\H$ such that  $\inf_{|x|_\H\le r} \E 1_K(X_{t_0}^x)>0.$
    \item[(iii)]  (In place of (2.4) in \cite{GM}): There exist   constants $\bb,k,c>0$ such that  $\E V(X_t^x)\le k V(x) \e^{-\bb t}+c,\ \ t\ge 0, x\in\H.$ \end{enumerate}

Firstly, according to \cite[Theorem 2.3.1]{Wbook} (see also \cite[Theorem 1.1]{LW} under a more specific framework), for any $p>1$ there exists a continuous function $\Psi_p$ on $\H\times\H\times (0,\infty)$ with $\Psi_p(x,x,t)=0$ such that
the Harnack inequality
\beq\label{H0} |P_t f(x)|^p\le (P_t|f|^p)(y) \e^{\Psi_p(x,y,t)},\ \ x,y\in \H, t>0, f\in \B_b(\H) \end{equation} holds. By \cite[Theorem 1.4.1]{Wbook} (see also \cite[Proposition 3.1]{WY}) for $P=P_t, \Psi=\Psi_p(\cdot,t)$ and $\Phi(s)=s^p$, this implies that $P_t$ has a unique   invariant probability measure $\mu$,   $P_t$ is strong Feller and has  a strictly positive density with respect to $\mu$.  Moreover, by the continuity of the solution, the Harnack inequality \eqref{H0} also implies that $\mu$ has full support on $\H$ (see the proof of Corollary 1.3(1) in \cite{WX}).  Therefore, (i) holds.

Next, since $h$ has relatively compact level sets in $\H$, it follows from \eqref{CM} that $P_{t_0}1_K(0)>0$ holds for some $t_0>0$ and compact set $K$ in $\H$. Indeed,  \eqref{CM} implies $c_0:=\E h(X_{t_0})<\infty$ for some $t_0>0$, so that we may take $K$ being the closure of $\{h\le c_0+1\}$.
Then it follows from \eqref{H0} that for any $r>0$,
$$\inf_{|x|_\H\le r} P_{t_0} 1_K(x) \ge (P_{t_0}1_K(0))^p \inf_{|x|_\H\le r} \e^{-\Phi_p(0,x,t_0)}>0.$$Thus, (ii) holds.

Finally, by \eqref{5.1} and the It\^o formula, we have
$$\d |X_t|_\H^2\le \big(2\aa+\|Q\|_{HS}^2 -2\eta |X_t|_\H^{1+r}\big)\d t +2\<X_t, Q\d W_t\>_\H.$$  Then for any $\gg >0$,
\beg{equation*}\beg{split} &\d \e^{\gg (1+|X_t|_\H^2)^{(1-r)/2}} \\
&\le  \ff{\gg(1-r)}2 \e^{\gg (1+|X_t|_\H^2)^{(1-r)/2}} (1+|X_t|_\H^2)^{-\ff{1+r}2} \d |X_t|_\H^2\\
 &\qquad + 2\gg^2(1+|X_t|_\H^2)^{-(1+r)} \|Q\|^2_{E\to \H} |X_t|_\H^2 \e^{\gg (1+|X_t|_\H^2)^{(1-r)/2}} \d t\\
&\le \gg \e^{\gg (1+|X_t|_\H^2)^{(1-r)/2}}\bigg\{\ff{(1-r)(2\aa+\|Q\|_{HS}^2)}{2(1+|X_T|_\H^2)^{\ff{1+r}2}} +\ff{2\gg \|Q\|_{E\to\H}^2}{(1+|X_t|_\H^2)^r}
-\ff{\eta (1-r)|X_t|_\H^{1+r}}{(1+|X_t|_\H^2)^{\ff{1+r}2}}\bigg\}\d t+\d M_t\\
&\le \big\{C_1-C_2 \e^{\gg (1+|X_t|_\H^2)^{(1-r)/2}}\big\}\d t +\d M_t\end{split}\end{equation*} for some constants $C_1,C_2>0$ and some local martingale $M_t$. This implies
$$\d \e^{C_2t+\gg (1+|X_t|_\H^2)^{(1-r)/2}}\le C_1 \e^{C_2 t}\d t + \e^{C_2 t}\d M_t.$$ Letting $\tau_n\uparrow\infty$ be a sequence of stopping times such that $(M_{t\land \tau_n})_{t\ge 0}$ is a martingale for every $n\ge 1$, we obtain
\beg{equation*}\beg{split} &\e^{C_2 t} \E V(X_t)=\e^{C_2 t} \E  \e^{\gg (1+|X_t|_\H^2)^{(1-r)/2}}= \E \liminf_{n\to \infty} \e^{C_2t\land \tau_n+ \gg (1+|X_{t\land\tau_n}|_\H^2)^{(1-r)/2}}\\
&\le \liminf_{n\to\infty}\E  \e^{C_2t\land \tau_n+ \gg (1+|X_{t\land\tau_n}|_\H^2)^{(1-r)/2}}
 \le \e^{\gg(1+|X_0|^2)^{(1-r)/2}} + \liminf_{n\to\infty}\E \int_0^{t\land\tau_n} C_1\e^{C_2 s}\d s\\
  &= \e^{\gg(1+|X_0|^2)^{(1-r)/2}} +\ff{C_1(\e^{C_2 t}-1)}{C_2}\le  V(X_0)+ \ff{C_1}{C_2}\e^{C_2 t},\ \ X_0\in\H.\end{split}\end{equation*}
 this implies (iii) for some $\bb=C_2, k=1$ and $c=\ff{C_1}{C_2}.$
\end{proof}

To illustrate Theorem \ref{T5.1}, we let $\DD,\m, \H:=H^{-1}, \{\ll_i,e_i\}_{i\ge 1}, E:=L^2(\m)$ and $W_t$ be in Section 3, and  consider the equation
 \beq\label{5E} \d X_t =\DD X_t^{r}\d t +Q\d W_t\end{equation} for  some $r\in (0,1)$, and $Qe_i:= q_i e_i (i\ge 1)$ with $\{q_i\}_{i\ge 1}\subset\R$ satisfying
 \beq\label{N}\|Q\|_{HS}^2:= \sum_{i=1}^\infty \ff{q_i^2}{\ll_i} <\infty,\ \ \inf_{i\ge 1} |q_i| \ll_i ^{\ff{1-\vv}{\theta}-\ff 1 2}> 0\end{equation} for some constants $\theta\ge \ff 4 {r+1} $ and $\vv\in (\ff{1-r}{2(1+r)}, 1).$ Since $\ll_i=\ff{\pi^2i^2}{l^2}$, if  $r\in (\ff 1 3,1)$ then for any $\kk\in (\ff 1 4, \ff{1+3r}8)\ne\emptyset$,   $q_i:= \ll_i^{\ff 1 2-\kk} (i\ge 1)$    satisfies \eqref{N} for  $\theta= \ff 4 {r+1} $ and $\vv=1-\ff{4\kk}{1+r}\in (\ff{1-r}{2(1+r)}, 1).$

 \beg{cor} Let $P_t$ be the Markov semigroup associated to $\eqref{5E}$ such that $\eqref{N}$ holds for some constants $\theta\ge \ff 4 {r+1} $ and $\vv\in (\ff{1-r}{2(1+r)}, 1)$. Then the assertion of Theorem $\ref{T5.1}$ holds. \end{cor}

 \beg{proof}
 By \eqref{N} we have $Q\in \L_{HS}(E,\H)$. Then it is easy to see that assumptions {\bf (A1)}-{\bf (A4)} hold for $b(u):= \DD u^r$ for $u\in \V:= L^{1+r}(\m)$, provided $\V$ is continuously embedded into $\H$. In fact, according to the proof of Corollary 3.2 in \cite{LW}, since $d:=1\in (0,\ff{2\vv(r+1)}{1-r})$  due to \eqref{N}, the classical Nash inequality
 $$\|f\|_2^{2+\ff 4 d}\le C\m(|\nn f|^2),\ \ f\in C_0^1((0,l)), \m(|f|)=1$$ for some constant $C>0$ implies that
 \beq\label{F1} \|x\|_{r+1}^2\ge c \sum_{i=1}^\infty \ff{\m(xe_i)^2}{\ll_i^\vv},\ \ x\in \V\end{equation}holds for some constant $c>0.$ Since $\vv<1$, this implies that $\V$ is   compactly (hence, also continuously) embedded into $\H$. So, it remains to verify conditions \eqref{5.1} and \eqref{5.2} in Theorem \ref{T5.1} for $h(u):= \|u\|_{r+1}$ and some constants $\aa,\eta>0.$

Since by \eqref{F1} we have $h(u)^2:=\|u\|_{r+1}^2\ge c \ll_1^{1-\vv} \|u\|_\H^2,$ \eqref{5.1} with   some $\eta>0$ and any $\aa>0$ follows from the fact that
$_{\V^*}\<b(u),u\>_\V= -h(u)^{1+r}.$ Next, by \eqref{N} we have
$|q_i|\ge c_1\ll_i^{\ff 1 2-\ff{1-\vv}\theta}$ for some constant $c_1>0$ and all $i\ge 1$. Combining this with \eqref{F1} we obtain \beg{equation}\label{F2}\beg{split}  \|x\|_Q^\theta &=\Big(\sum_{i\ge 1}\ff{\m(x e_i)^2}{q_i^2}\Big)^{\ff\theta 2}
  \le \Big(\sum_{i\ge 1} \ff{\ll_i^{\ff{\theta-2}2}\m(xe_i)^2}{|q_i|^\theta}\Big)\Big(\sum_{i\ge 1}\ff{\m(xe_i)^2}{\ll_i}\Big)^{\ff{\theta-2}2}\\
  &\le c_2|x|_\H^{\theta-2} \sum_{i\ge 1} \ff{\m(x e_i)^2}{\ll_i^\vv}\le c_3|x|_\H^{\theta-2} \|x\|_{r+1}^2,\ \ x\in\V\end{split}\end{equation} for some constants $c_2,c_3>0$. Moreover, by the H\"older inequality and noting that
  $$\m\big((|u|\lor |v|)^{1+r}\big)\le h(u)^{1+r}+h(v)^{1+r}\le 2 \{h(u)\lor h(v)\}^{1+r},$$ we obtain
  \beg{equation*}\beg{split} \|u-v\|_{r+1}^{1+r} &:=\m(|u-v|^{1+r}) \le \m\big(|u-v|^2(|u|\lor |v|)^{r-1}\big)^{\ff{1+r}2} \m\big((|u|\lor |v|)^{1+r}\big)^{\ff{1-r}2}\\
  &\le 2^{\ff{1-r}2}   \m\big(|u-v|^2(|u|\lor |v|)^{r-1}\big)^{\ff{1+r}2}\{h(u)\lor h(v)\}^{\ff{1-r^2}2}.\end{split}\end{equation*}
Combining this with \eqref{F2} we arrive at
\beg{equation*}\beg{split}& _{\V^*}\<b(u)-b(v),u-v\>_\V=-\m\big((u^r-v^r)(u-v)\big)\le -r \m\big(|u-v|^2(|u|\lor |v|)^{r-1}\big)\\
&\le -\ff{r 2^{\ff{r-1}2}\|u-v\|_{r+1}^2}{(h(u)\lor h(v))^{1-r}}\le\ff{-\eta \|u-v\|_Q^\theta}{|u-v|_\H^{\theta-2}(h(u)\lor h(v))^{1-r}}\end{split}\end{equation*} for some constant $\eta>0$ and all $u,v\in\V$. Thus, \eqref{5.2} holds.
  \end{proof}
  
\paragraph{Acknowledgement.} The author would like to thank the referees for their corrections and useful comments. 

\beg{thebibliography}{99}

\bibitem{AP} D. G. Aronson, L. A. Peletier, \emph{Large time behaviour of solutions of the porous medium equation in bounded domains,} J. Diff. Equ. 39(1981), 378--412.
\bibitem{CG} A. Chojnowska-Michalik, B. Goldys, \emph{Nonsymmetric Ornstein-Uhlenbeck semigroup as second quantized operator,} J. Math. Kyoto Univ. 36(1996), 481--498.

\bibitem{DRRW} G. Da Prato, M. R\"ockner, B.L. Rozovskii, F.-Y. Wang, \emph{Strong solutions of Generalized porous media equations: existence, uniqueness and ergodicity, } Comm. Part. Diff. Equ. 31 (2006), 277--291.


\bibitem{GM} B. Goldys, B. Maslowski, \emph{Exponential ergordicity for stochastic reaction-diffusion equations,} in $``$Stochastic Partial Differential Equations and Applications VII." Lecture Notes Pure Appl. Math. 245(2004), 115--131. Chapman Hall/CRC Press.

\bibitem{GM2} B. Goldys, B. Maslowski, \emph{Lower estimates of transition density and bounds on exponential ergodicity for stochastic PDEs,} Ann. Probab. 34(2006), 1451--1496.


\bibitem{KR}
  N.V. Krylov, B.L. Rozovskii,
  \emph{Stochastic evolution equations,}
  Translated from Itogi Naukii Tekhniki, Seriya Sovremennye Problemy Matematiki,
  14(1979), 71--146, Plenum Publishing Corp. 1981.


\bibitem{L} W. Liu, \emph{Harnack inequality and applications for stochastic evolution equations with monotone drifts,}  J. Evol. Equ. 9(2009),  747--770.
    
\bibitem{Liu11} W. Liu, \emph{Ergodicity of transition semigroups for stochastic fast diffusion equations,} Front. Math. China 6(2011), 449--472.

\bibitem{LR} W. Liu, M. R\"ockner, \emph{SPDE in Hilbert space with locally monotone coefficients,} J. Funct. Anal.   259(2010), 2902--2922.

\bibitem{LR2} W. Liu, M. R\"ockner, \emph{Local and global well-posedness of SPDE with generalized coercivity conditions,} J. Differential Equations  254(2013), 725--755.
    
\bibitem{LT11} W. Liu, J.M. T\"olle, \emph{Existence and uniqueness of invariant measures for stochastic evolution equations with weakly dissipative drifts,} Elect. Comm. Probab. 16(2011), 447--457. 
    
\bibitem{LW} W. Liu, F.-Y. Wang, \emph{Harnack inequality and strong Feller property for stochastic fast-diffusion equations,} J. Math. Anal. Appl. 342(2008), 651--662.

\bibitem{P1} E. Pardoux, \emph{Sur des equations aux d\'eriv\'ees
partielles stochastiques monotones,} C. R. Acad. Sci. 275(1972),
A101--A103.

\bibitem{P2} E. Pardoux, \emph{Equations aux d\'eriv\'ees
partielles stochastiques non lineaires monotones: Etude de solutions
fortes de type Ito,} Th\'ese Doct. Sci. Math. Univ. Paris Sud. 1975.

\bibitem{RRW} J. Ren, M. R\"ockner, F.-Y. Wang, \emph{Stochastic generalized porous media and fast diffusion equations,} J. Differential Equations 238(2007), 118--152.

\bibitem{RW} M. R\"ockner, F.-Y. Wang, \emph{Harnack and functional inequalities for generalized Mehler semigroups,}  J. Funct. Anal.  203(2003), 237--261.
\bibitem{RW08}  M. R\"ockner, F.-Y. Wang, \emph{ Non-monotone stochastic generalized porous media equations,} J. Differential Equations  245(2008), 3898-3935.
\bibitem{W07} 	F.-Y. Wang, \emph{Harnack inequality and applications for stochastic generalized porous media equations,}  Annals of Probability 35(2007), 1333--1350.

\bibitem{Wbook} F.-Y. Wang, \emph{Harnack Inequality and Applications for Stochastic Partial Differential Equations,} Springer Briefs in Mathematics,  Springer, New York, 2013.
    
\bibitem{WX} F.-Y. Wang,   L. Xu, \emph{Derivative formula and applications for hyperdissipative stochastic Navier-Stokes/Burgers equations,} Infin. Dimens. Anal. Quantum Probab. Relat. Top. 15(2012), no. 3, 1250020, 19 pp.

\bibitem{WY} F.-Y. Wang, C. Yuan, \emph{Harnack inequalities for functional SDEs with multiplicative noise and applications,} Stoch. Proc. Appl. 121(2011), 2692--2710.

\end{thebibliography}

\end{document}